\input amstex
\input amsppt.sty
\hsize=12.5cm
\vsize=540pt
\baselineskip=16truept
\NoRunningHeads
\TagsOnRight
 \define\CC{\Bbb C}
 \define\CF{\Cal F}
 \define\CO{\Cal O}
 \redefine\Delta{\varDelta}
 \define\eps{\varepsilon}
 \define\halfskip{\vskip6pt plus 1pt}
 \define\intt{\operatorname{int}}
 \font\msbm=msbm10
 \define\MX{\text{\msbm X}}
 \define\MT{\text{\msbm T}}
 \define\NN{\Bbb N}
 \redefine\phi{\varphi}
 \define\PSH{\Cal P\Cal S\Cal H}

 \define\tdots{\times\dots\times}
 \define\too{\longrightarrow}
 \define\wdht{\widehat}
 \define\wdtl{\widetilde}
\catcode`@=11
\def\qed{\ifhmode\textqed\else\ifinner\quad\square
   \else\eqno\square\fi\fi}
\def\textqed{{\unskip\nobreak\penalty50
    \quad\hbox{}\nobreak\hfil$\m@th\square$
    \parfillskip=0pt \finalhyphendemerits=0\par}}
\define\ftnt#1{\global\advance\footmarkcount@1%
\footnote"$\;\left({}^{\number\footmarkcount@}\right)$"{\; #1}}
\catcode`@=\active

\document
\topmatter
\title
An extension theorem for separately holomorphic functions with pluripolar
singularities
\endtitle

\author
Marek Jarnicki (Krak\'ow)\;
\footnote"${}^{(\dag)}$"
{\;\hbox{Research partially supported by the KBN grant No\. 5 P03A 033 21.}}
Peter Pflug (Oldenburg)\;
\footnote"${}^{(\ddag)}$"
{\;Research partially supported by the Nieders\"achsisches
Ministerium f\"ur Wissenschaft und
\hbox{\hskip12pt Kultur, Az. 15.3 -- 50 113(55) PL.}}
\endauthor

\address
\vbox{ \hbox{} \hbox{Uniwersytet Jagiello\'nski} \hbox{Instytut
Matematyki} \hbox{30-059 Krak\'ow, Reymonta 4,} \hbox{Poland}
\hbox{{\it E-mail: }{\rm jarnicki{\@}im.uj.edu.pl}} } \hfill
\vbox{ \hbox{Carl von Ossietzky Universit\"at Oldenburg}
\hbox{Fachbereich Mathematik} \hbox{Postfach 2503} \hbox{D-26111
Oldenburg, Germany} \hbox{{\it E-mail: }{\rm
pflug{\@}mathematik.uni{-oldenburg.de}}} }
\endaddress

\abstract
Let $D_j\subset\CC^{n_j}$ be a pseudoconvex domain and let $A_j\subset D_j$
be a locally pluriregular set, $j=1,\dots,N$. Put
$$
X:=\bigcup_{j=1}^N A_1\tdots A_{j-1}\times D_j\times A_{j+1}\tdots A_N
\subset\CC^{n_1}\tdots\CC^{n_N}=\CC^n.
$$
Let $U\subset\CC^n$ be an open neighborhood of $X$ and let $M\subset U$ be a
relatively closed subset of $U$. For $j\in\{1,\dots,N\}$ let
$\Sigma_j$ be the set of all $(z',z'')\in(A_1\tdots A_{j-1})
\times(A_{j+1}\tdots A_N)$ for which the fiber
$M_{(z',\cdot,z'')}:=\{z_j\in\CC^{n_j}\: (z',z_j,z'')\in M\}$ is not pluripolar.
Assume that $\Sigma_1,\dots,\Sigma_N$ are pluripolar. Put
$$
X':=\bigcup_{j=1}^N\{(z',z_j,z'')\in(A_1\tdots A_{j-1})\times D_j
\times(A_{j+1}\tdots A_N)\: (z',z'')\notin\Sigma_j\}.
$$
Then there exists a relatively closed pluripolar
subset $\wdht M\subset\wdht X$ of the `envelope of holomorphy' $\wdht X\subset\CC^n$
of $X$ such that:

$\bullet$ $\wdht M\cap X'\subset M$,

$\bullet$ for every function $f$ separately holomorphic on $X\setminus M$
there exists exactly one function $\wdht f$ holomorphic on
$\wdht X\setminus\wdht M$ with $\wdht f=f$ on $X'\setminus M$, and

$\bullet$ $\wdht M$ is singular with respect to the family of all functions
$\wdht f$.

\noindent Some special cases were previously studied in \cite{Jar-Pfl 2001c}.
\endabstract

\subjclass 32D15, 32D10\endsubjclass

\endtopmatter


\noindent{\bf 1. Introduction. Main Theorem.} Let $N\in\NN$, $N\geq2$, and let
$$
\varnothing\neq A_j\subset D_j\subset\CC^{n_j},
$$
where $D_j$ is a domain, $j=1,\dots,N$. We define an {\it $N$--fold cross}
$$
\align
X&=\MX(A_1,\dots,A_N;D_1,\dots,D_N)\\
&:=\bigcup_{j=1}^NA_1\tdots A_{j-1}\times D_j\times A_{j+1}\tdots A_N
\subset\CC^{n_1+\dots+n_N}=\CC^n.
\endalign
$$
Observe that $X$ is connected.

\halfskip

Let $\Omega\subset\CC^n$ be an open set and let $A\subset\Omega$. Put
$$
h_{A,\Omega}:=\sup\{u\:\; u\in\PSH(\Omega),\;u\leq1 \text{ on }
\Omega,\; u\leq0 \text{ on } A\},
$$
where $\PSH(\Omega)$ denotes the set of all functions plurisubharmonic on
$\Omega$. Define
$$
\omega_{A,\Omega}:=\lim_{k\to+\infty}h^\ast_{A\cap\Omega_k,\Omega_k},
$$
where $(\Omega_k)_{k=1}^\infty$ is a sequence of relatively compact open sets
$\Omega_k\subset\Omega_{k+1}\subset\subset\Omega$ with
$\bigcup_{k=1}^\infty\Omega_k=\Omega$ ($h^\ast$ denotes the upper
semicontinuous regularization of $h$). Observe that the definition is
independent of the exhausting sequence $(\Omega_k)_{k=1}^\infty$.
Moreover, $\omega_{A,\Omega}\in\PSH(\Omega)$. Recall that if $\Omega$ is
bounded, then $\omega_{A,\Omega}=h^\ast_{A,\Omega}$.

\halfskip

For an $N$--fold cross $X=\MX(A_1,\dots,A_N;D_1,\dots,D_N)$ put
$$
\wdht X:=\{(z_1,\dots,z_N)\in D_1\tdots D_N\:
\sum_{j=1}^N\omega_{A_j,D_j}(z_j)<1\}.
$$
Observe that if $D_1,\dots,D_N$ are pseudoconvex, then $\wdht X$ is a
pseudoconvex open set in $\CC^n$.

\halfskip

We say that a subset $\varnothing\neq A\subset\CC^n$ is {\it locally
pluriregular} if $h^\ast_{A\cap\Omega,\Omega}(a)=0$ for any $a\in A$ and for
any open neighborhood $\Omega$ of $a$ (in particular, $A\cap\Omega$ is
non-pluripolar).

Note that if $A_1,\dots,A_N$ are locally pluriregular, then $X\subset\wdht X$
and $\wdht X$ is connected (\cite{Jar-Pfl 2001c}, Lemma 4).

\halfskip

Let $U$ be an open neighborhood of $X$ and let $M\subset U$ be a relatively
closed set. We say that a function $f\:X\setminus M\too\CC$ is {\it separately
holomorphic} ($f\in\CO_s(X\setminus M)$) if for any
$(a_1,\dots,a_N)\in A_1\tdots A_N$ and $j\in\{1,\dots,N\}$ the function
$f(a_1,\dots,a_{j-1},\cdot,a_{j+1},\dots,a_N)$ is holomorphic in the open set
$$
D_j\setminus M_{(a_1,\dots,a_{j-1},\cdot,a_{j+1},\dots,a_N)},
$$
where
$$M_{(a_1,\dots,a_{j-1},\cdot,a_{j+1},\dots,a_N)}:=
\{z_j\in\CC^{n_j}\: (a_1,\dots,a_{j-1},z_j,a_{j+1},\dots,a_N)\in M\}.
$$

\halfskip

Suppose that $S_j\subset A_1\tdots A_{j-1}\times A_{j+1}\tdots A_N$,
$j=1,\dots,N$, and define the {\it generalized $N$--fold cross}
$$
\multline
T=\MT(A_1,\dots,A_N;D_1,\dots,D_N;S_1,\dots,S_N)\\
:=\bigcup_{j=1}^N\{(z',z_j,z'')\in(A_1\tdots A_{j-1})\times D_j
\times(A_{j+1}\tdots A_N)\: (z',z'')\notin S_j\}.
\endmultline
$$
It is clear that $T\subset X$. Observe that
$$
\MX(A_1,\dots,A_N;D_1,\dots,D_N)=
\MT(A_1,\dots,A_N;D_1,\dots,D_N;\varnothing,\dots,\varnothing).
$$
Moreover, if $N=2$, then $\MT(A_1,A_2;D_1,D_2;S_1,S_2)=\MX(A_1\setminus S_2,
A_2\setminus S_1;D_1,D_2)$. Consequently, any generalized $2$--fold cross is
a $2$--fold cross.

\halfskip

Let $S\subset\Omega$ be a relatively closed pluripolar subset of
an open set $\Omega\subset\CC^n$. Let
$\CF\subset\CO(\Omega\setminus S)$. We say that $S$ is {\it
singular with respect to $\CF$} if for any point $a\in S$ there
exists a function $f_a\in\CF$ which is not holomorphically
extendible to a neighborhood of $a$ (cf\. \cite{Jar-Pfl 2000},
\S\;3.4). Equivalently: the set $S$ is minimal in the sense that
there is no relatively closed set $S'\varsubsetneq S$ such that
any function from $\CF$ extends holomorphically to
$\Omega\setminus S'$. It is clear that for any relatively closed
pluripolar set $S\subset\Omega$ and for any family
$\CF\subset\CO(\Omega\setminus S)$ there exists a relatively
closed set $S'\subset S$ such that any function $f\in\CF$ extends
to an $f'\in\CO(\Omega\setminus S')$ and $S'$ is singular with
respect to the family $\{f'\:f\in\CF\}$.

\halfskip

The main result of our paper is the following extension theorem for
separately holomorphic functions.

\proclaim{Main Theorem} Let $D_j\subset\CC^{n_j}$ be a pseudoconvex domain,
let $A_j\subset D_j$ be a locally pluriregular set, $j=1,\dots,N$, and let
$U$ be an open neighborhood of the $N$--fold cross
$$
X:=\MX(A_1,\dots,A_N;D_1,\dots,D_N).
$$
Let $M\subset U$ be a relatively closed subset of $U$ such that
for each $j\in\{1,\dots,N\}$ the set
$$
\multline
\Sigma_j=\Sigma_j(A_1,\dots,A_N;M)\\
:=\{(z',z'')\in(A_1\tdots A_{j-1})\times(A_{j+1}\tdots A_N)\:
M_{(z',\cdot,z'')} \text{ is not pluripolar}\}
\endmultline
$$
is pluripolar. Put
$$
X':=\MT(A_1,\dots,A_N;D_1,\dots,D_N;\Sigma_1,\dots,\Sigma_N).
$$
Then there exists a relatively closed pluripolar
set $\wdht M\subset\wdht X$ such that:

$\bullet$ $\wdht M\cap X'\subset M$,

$\bullet$ for every $f\in\CO_s(X\setminus M)$ there exists exactly one
$\wdht f\in\CO(\wdht X\setminus\wdht M)$ with $\wdht f=f$ on
$X'\setminus M$,

$\bullet$ $\wdht M$ is singular with respect to the family
$\{\wdht f\: f\in\CO_s(X\setminus M)\}$, and

$\bullet$ $\wdht X\setminus\wdht M$ is pseudoconvex.

In particular, $\wdht X\setminus\wdht M$ is the envelope of holomorphy of
$X\setminus M$ with respect to the space of separately holomorphic functions.
\endproclaim

Notice that if $M\subset U$ is a pluripolar set,
then $\Sigma_1,\dots,\Sigma_N$ are always pluripolar (cf\. Lemma 8(a)).

The case where $N=2$, $n_1=n_2=1$, $D_1=D_2=\CC$ was studied in
\cite{Jar-Pfl 2001b}, Th.~2.

\proclaim{Corollary 1} Let $D_j$, $A_j$, $j=1,\dots,N$, $X$, and $U$ be as
in the Main Theorem. Assume that $M\subset U$ is a relatively closed set such that
for any $(a_1,\dots,a_N)\in A_1\tdots A_N$ and $j\in\{1,\dots,N\}$ the fiber
$M_{(a_1,\dots,a_{j-1},\cdot,a_{j+1},\dots,a_N)}$
is pluripolar\; \ftnt{That is $\Sigma_1=\dots=\Sigma_N=\varnothing$.}.
Then there exists a relatively closed pluripolar
set $\wdht M\subset\wdht X$ such that:

$\bullet$ $\wdht M\cap X\subset M$,

$\bullet$ for every $f\in\CO_s(X\setminus M)$ there exists exactly one
$\wdht f\in\CO(\wdht X\setminus\wdht M)$ with $\wdht f=f$ on
$X\setminus M$, and

$\bullet$ the domain $\wdht X\setminus\wdht M$ is pseudoconvex.
\endproclaim

The case where $N=2$, $D_2=\CC^{n_2}$, and $A_2$ is open was studied in
\cite{Chi-Sad 1988} (for $n_2=1$) and \cite{Kaz 1988} (for arbitrary $n_2$).

The proof of the Main Theorem will be presented in Sections 3 (for $N=2$) and
4 (for arbitrary $N$).

\halfskip

The following two examples illustrate the role played by the sets
$\Sigma_j$ and show that the assertion of the Main Theorem is in some sense
optimal.

\proclaim{Example 2}\rm Let $n_1=n_2=1$, $D_1=D_2=\CC$, $A_1=E:=$
the unit disc.

(a) Let $A_2:=E$, $X:=\MX(E,E;\CC,\CC)=(E\times\CC)\cup(\CC\times
E)$, and $M:=\{0\}\times\overline E$. Then $\Sigma_1=\varnothing$,
$\Sigma_2=\{0\}$, $X'=\MX(E\setminus\{0\},E;\CC,\CC)$, $\wdht
M=\{0\}\times\CC$.

Put $f_0(z,w):=1/z$, $z\neq0$, and $f_0(0,w)=1$, $|w|>1$. Then
$f_0\in\CO_s(X\setminus M)$ and $\wdht M$ is singular with respect
to $f_0$.

(b) Let $A_2:=E\setminus r\overline E$, $X:=\MX(E,E\setminus r\overline
E;\CC,\CC)$, and $M:=\{0\}\times\{|w|=r\}$ for some $0<r<1$. Then
$\Sigma_1=\varnothing$, $\Sigma_2=\{0\}$,
$X'=\MX(E\setminus\{0\},A_2;\CC,\CC)$, $\wdht M=\varnothing$.

Put
$$
f_0(z,w):=\cases w &\text{if $z\neq0$ or ($z=0$ and $|w|>r$)}\\
                 0 &\text{if \hskip43pt $z=0$ and $|w|<r$}
\endcases,\quad (z,w)\in X\setminus M.
$$
Then $f_0\in\CO_s(X\setminus M)$, $\wdht f_0(z,w)\equiv w$, and
$\wdht f_0(0,w)\neq f(0,w)$, $0<|w|<r$.
\endproclaim

\halfskip

\noindent{\bf 2. Auxiliary results.}

In the case $M=\varnothing$ the problem of extension of separately holomorphic
functions was studied by many authors (under various assumptions on
$(D_j,A_j)_{j=1}^N$), e.g\. \cite{Sic 1969}, \cite{Zah 1976}, \cite{Sic 1981},
\cite{Shi 1989}, \cite{Ngu-Zer 1991}, \cite{Ngu 1997}, \cite{Ale-Zer 2001}
(for $N=2$), and \cite{Sic 1981}, \cite{Ngu-Zer 1995}, \cite{Jar-Pfl 2001c}
(for arbitrary $N$).

\proclaim\nofrills{Theorem 3}\; {\rm (\cite{Ngu-Zer 1995},
\cite{Ale-Zer 2001}).}
Let $(D_j,A_j)_{j=1}^N$ and $X$ be as in the Main Theorem. Then any function
from $\CO_s(X)$ extends holomorphically to the pseudoconvex domain $\wdht X$.
\endproclaim

\halfskip

The case where $M$ is analytic was studied in \cite{\"Okt 1998},
\cite{\"Okt 1999}, \cite{Sic 2000}, \cite{Jar-Pfl 2001a}. The problem was
completely solved in \cite{Jar-Pfl 2001c}.

\proclaim\nofrills{Theorem 4}\;{\rm (\cite{Jar-Pfl 2001b}).}
Let $(D_j,A_j)_{j=1}^N$ and $X$ be as in the Main Theorem. Let
$M\varsubsetneq U$ be an analytic subset of an open connected neighborhood
$U$ of $X$. Then there exists an analytic set $\wdht M\subset\wdht X$ such
that:

$\bullet$ $\wdht M\cap U_0\subset M$ for an open neighborhood $U_0$ of $X$,
$U_0\subset U$,

$\bullet$ for every $f\in\CO_s(X\setminus M)$ there exists exactly one
$\wdht f\in\CO(\wdht X\setminus\wdht M)$ with $\wdht f=f$ on $X\setminus M$,
and

$\bullet$ the domain $\wdht X\setminus\wdht M$ is pseudoconvex.
\endproclaim

\proclaim{Remark 5}\rm It is an natural idea to try to obtain Theorem 4 from
the Main Theorem. More precisely, let $(D_j,A_j)_{j=1}^N$, $X$, $U$, and $M$ be
as in Theorem 4. Then, by the Main Theorem, there exists a relatively closed
pluripolar set $\wdht M\subset\wdht X$ which has all the properties listed
in the Main Theorem.
We would like to know whether there is a direct argument showing that
$\wdht M$ must be analytic.
\endproclaim

\halfskip

The following two results will play the fundamental role in the sequel.

\proclaim\nofrills{Theorem 6}\;{\rm (\cite{Chi 1993}).}
Let $D\subset\CC^n$ be a domain and let $\wdht D$ be the envelope of
holomorphy of $D$. Assume that $S$ is relatively closed pluripolar subset of
$D$. Then there exists a relatively closed pluripolar subset $\wdht S$
of $\wdht D$ such that $\wdht S\cap D\subset S$ and
$\wdht D\setminus\wdht S$ is the envelope of holomorphy of $D\setminus S$.
\endproclaim

\halfskip

\proclaim\nofrills{Theorem 7}\;{\rm (\cite{Jar-Pfl 2001b}).}
Let $A\subset E^{n-1}$ be locally pluriregular, let
$$
X:=\MX(A,E;E^{n-1},\CC)
$$
(notice that $\wdht X=E^{n-1}\times\CC$), and let $U\subset E^{N-1}\times\CC$ be an open
neighborhood of $X$. Let $M\subset U$ be a relatively closed set such that
$M\cap E^n=\varnothing$ and for any $a\in A$ the fiber $M_{(a,\cdot)}$ is
polar. Then there exists a relatively closed pluripolar set
$S\subset E^{n-1}\times\CC$ such that

$\bullet$ $S\cap X\subset M$,

$\bullet$ any function from $\CO_s(X\setminus M)$ extends
holomorphically to $E^{n-1}\times\CC\setminus S$, and

$\bullet$ $E^{n-1}\times\CC\setminus S$\; \ftnt{Here
and in the sequel to simplify notation we write $P_1\tdots P_k\setminus Q$
instead of $(P_1\tdots P_k)\setminus Q$.} is pseudoconvex.
\endproclaim

Notice that the above result is a special case of our Main Theorem with
$N=2$, $n_1=n-1$, $D_1=E^{n-1}$, $A_1=A$, $n_2=1$, $D_2=\CC$, $A_2=E$,
$\Sigma_1=\Sigma_2=\varnothing$.

\demo{Proof} It is known (cf\. \cite{Chi-Sad 1988}) that each function
$f\in\CO_s(X\setminus M)$ has the univalent domain of existence
$G_f\subset E^{n-1}\times\CC$\;
\ftnt{We like to thank Professor Evgeni Chirka for
explaining us some details of the proof of Theorem 1 in \cite{Chi-Sad 1988}.}.
Let $G$ denote the connected component of
$\intt\bigcap_{f\in\CO_s(X\setminus M)}G_f$ that
contains $E^n$ and let $S:=E^{n-1}\times\CC\setminus G$. It remains to show
that $S$ is pluripolar.

Take $(a,b)\in A\times\CC\setminus M$. Since $M_{(a,\cdot)}$
is polar, there exists a curve $\gamma\:[0,1]\too\CC\setminus M_{(a,\cdot)}$
such that $\gamma(0)=0$, $\gamma(1)=b$. Take an $\eps>0$ so small that
$$
\Delta_a(\eps)\times(\gamma([0,1])+\Delta_0(\eps))\subset U\setminus M,
$$
where $\Delta_{z_0}(r)=\Delta_{z_0}^k(r)\subset\CC^k$ denotes the polydisc
with center $z_0\in\CC^k$ and radius $r>0$.
Put $V_b:=E\cup(\gamma([0,1])+\Delta_0(\eps))$ and consider the cross
$$
Y:=\MX(A\cap\Delta_a(\eps),E;\Delta_a(\eps),V_b).
$$
Then $f\in\CO_s(Y)$ for any $f\in\CO_s(X\setminus M)$. Consequently, by
Theorem 3, we get $\wdht Y\subset G_f$, $f\in\CO_s(X\setminus M)$.
Hence $\wdht Y\subset G$. In particular, we conclude that
$\{a\}\times(\CC\setminus M_{(a,\cdot)})\subset G$.

Thus $S_{(a,\cdot)}\subset M_{(a,\cdot)}$ for all $a\in A$.
Consequently, by Lemma 5 from \cite{Chi-Sad 1988}, $S$ is pluripolar.
\qed\enddemo

\proclaim{Lemma 8} {\rm (a)} Let $S\subset\CC^p\times\CC^q$ be pluripolar.
Then the set
$$
A:=\{z\in\CC^p\: S_{(z,\cdot)} \text{ is not pluripolar}\}
$$
is pluripolar.

\noindent{\rm (b)} Let $M\subset\CC^p\times\CC^q$ be such that for each
$a\in\CC^p$ the fiber $M_{(a,\cdot)}$ is pluripolar. Let
$C\subset\CC^p\times\CC^q$ be such that the set
$\{z\in\CC^p\: C_{(z,\cdot)} \text{ is not pluripolar}\}$ is not pluripolar
(e.g\. $C=C'\times C''$, where $C'\subset\CC^p$, $C''\subset\CC^q$ are
non-pluripolar). Then $C\setminus M$ is non-pluripolar.

\noindent{\rm (c)} Let $M\subset\CC^p\times\CC^q$ be such that for
each $a\in\CC^p$ the fiber $M_{(a,\cdot)}$ is pluripolar. Let
$A\subset\CC^p$ be locally pluriregular. Let
$C:=\{(a,b')\in A\times\CC^{q-1}\: M_{(a,b',\cdot)} \text{ is polar}\}$
Then $C$ is locally pluriregular.
\endproclaim

\demo{Proof} (a) Let $v\in\PSH(\CC^{p+q})$, $v\not\equiv-\infty$, be such that
$S\subset v^{-1}(-\infty)$. Define
$$
u(z):=\sup\{v(z,w)\:w\in\overline E^q\},\quad z\in\CC^p.
$$
Then $A\subset u^{-1}(-\infty)$. Moreover,
$u\in\PSH(\CC^p)$ and $u\not\equiv-\infty$.

(b) Suppose that $C\setminus M$ is pluripolar. Then, by (a), there exists a
pluripolar set $A\subset\CC^p$ such that the fiber
$(C\setminus M)_{(a,\cdot)}$ is pluripolar, $a\in\CC^p\setminus A$.
Consequently, the fiber $C_{(a,\cdot)}$ is pluripolar, $a\in\CC^p\setminus A$;
contradiction.

(c) Fix a point $(a_0,b'_0)\in C$ and a neighborhood
$U:=\Delta_{(a_0,b'_0)}(r)$. We have to show that
$h^\ast_{C\cap U,U}(a_0,b'_0)=0$. First we show that
$$
h^\ast_{C\cap U,U}(a_0,b'_0)\leq
h^\ast_{(A\cap\Delta_{a_0}(r))\times\Delta_{b'_0}(r),U}(a_0,b'_0).\tag*
$$
Indeed, let
$u\in\PSH(U)$ be such that $u\leq1$ and $u\leq0$ on $C\cap U$. Then for any
$a\in A\cap\Delta_{a_0}(r)$ the function $u(a,\cdot)$ is plurisubharmonic on
$\Delta_{b'_0}(r)$  and $u(a,\cdot)\leq 0$ on the set
$$
(C\cap U)_{(a,\cdot)}=\{b'\in\Delta_{b'_0}(r)\:
(M_{(a,\cdot)})_{(b',\cdot)} \text{ is polar}\}.
$$
By (a) (applied to the set $M_{(a,\cdot)}$), the set $\Delta_{b'_0}(r)
\setminus(C\cap U)_{(a,\cdot)}$ is pluripolar. Hence
$u(a,\cdot)\leq0$ on $\Delta_{b'_0}(r)$. Consequently,
$u\leq0$ on $(A\cap\Delta_{a_0}(r))\times\Delta_{b'_0}(r)$, which implies
that $h_{C\cap U,U}\leq h_{(A\cap\Delta_{a_0}(r))\times\Delta_{b'_0}(r),U}$,
and finally, $h^\ast_{C\cap U,U}(a_0,b'_0)\leq
h^\ast_{(A\cap\Delta_{a_0}(r))\times\Delta_{b'_0},U}(a_0,b'_0)$.

Now, in virtue of the product property of the relative extremal function
(cf\. \cite{Ngu-Sic 1991}), using (*) and the fact that $A$ is locally
pluriregular, we get
$$
\multline
h^\ast_{C\cap U,U}(a_0,b'_0)\leq
h^\ast_{(A\cap\Delta_{a_0}(r))\times\Delta_{b'_0}(r),U}(a_0,b'_0)=\\
\max\Big\{h^\ast_{A\cap\Delta_{a_0}(r),\Delta_{a_0}(r)}(a_0),\;
h^\ast_{\Delta_{b'_0}(r),\Delta_{b'_0}(r)}(b'_0)\Big\}
=h^\ast_{A\cap\Delta_{a_0}(r),\Delta_{a_0}(r)}(a_0)=0.
\endmultline
$$
\qed\enddemo

\proclaim{Lemma 9} Let $D_j$, $A_j$, $j=1,\dots,N$, and $X$ be as in the
Main Theorem. Let
$$
S_j\subset A_1\tdots A_{j-1}\times A_{j+1}\tdots A_N
$$
be pluripolar, $j=1,\dots,N$. Put
$$
T:=\MT(A_1,\dots,A_N;D_1,\dots,D_N;S_1,\dots,S_N).
$$
Then any function $f\in\CO_s(T)\cap\Cal C(T)$\;
\ftnt{We say that a function $f\:T\too\CC$ is {\it separately holomorphic}
if for any $j\in\{1,\dots,N\}$ and
$(a',a'')\in (A_1\tdots A_{j-1})\times(A_{j+1}\tdots A_N)\setminus S_j$ the
function $f(a',\cdot,a'')$ is holomorphic in $D_j$.}
extends holomorphically to $\wdht X$.
\endproclaim

If $N=2$, then the result is true for any function $f\in\CO_s(T)$
(see the proof).
In the case where $N\geq3$ we do not know whether the result is true for
arbitrary $f\in\CO_s(T)$.

\demo{Proof}
We apply induction on $N$. The case $N=2$ follows from Theorem 3 and the fact
that $\wdht X=\wdht T$ (recall that if $N=2$, then $T$ is a $2$--fold cross).
Moreover, if $N=2$, then the result is true for any $f\in\CO_s(T)$.

Assume that the result is true for $N-1\geq2$. Take an
$f\in\CO_s(T)\cap\Cal C(T)$. Let $Q$ denote the set of all $z_N\in A_N$ for
which there exists a $j\in\{1,\dots,N-1\}$ such that the fiber
$(S_j)_{(\cdot,z_N)}$ is not pluripolar. Then, by Lemma 8(a), $Q$
is pluripolar. Take a $z_N\in A_N\setminus Q$ and define
$$
T_{z_N}:=\MT(A_1,\dots,A_{N-1};D_1,\dots,D_{N-1};
(S_1)_{(\cdot,z_N)},\dots,(S_{N-1})_{(\cdot,z_N)}).
$$
Then $f(\cdot,z_N)\in\CO_s(T_{z_N})\cap\Cal C(T_{z_N})$. By the inductive
assumption, the function $f(\cdot,z_N)$ extends to an
$\wdht f_{z_N}\in\CO(\wdht Y)$, where
$Y=\MX(A_1,\dots,A_{N-1};D_1,\dots,D_{N-1})$.

Let $A':=A_1\times\dots\times A_{N-1}$. Consider the $2$--fold cross
$$
Z:=\MT(A',A_N;\wdht Y,D_N;S_N,Q)=((A'\setminus S_N)\times D_N)\cup
(\wdht Y\times(A_N\setminus Q)).
$$
Let $g\:Z\too\CC$ be given by the formulae:

$g(z',z_N):=f(z',z_N)$, $(z',z_N)\in(A'\setminus S_N)\times D_N$,

$g(z',z_N):=\wdht f_{z_N}(z')$, $(z',z_N)\in\wdht
Y\times(A_N\setminus Q)$.

Observe that $g$ is well-defined.

Indeed, let $(z',z_N)\in((A'\setminus S_N)\times D_N)\cap
(\wdht Y\times(A_N\setminus Q))$. If $z'\in T_{z_N}$, then obviously
$\wdht f_{z_N}(z')=f(z',z_N)$. Suppose that $z'\notin T_{z_N}$. Then
$$
\multline
z'\in P_{z_N}:=\bigcap_{j=1}^{N-1}\{(w',w_j,w'')\in(A_1\tdots A_{j-1})
\times A_j\times(A_{j+1}\tdots A_{N-1})\:\\
(w',w'')\in(S_j)_{(\cdot,z_N)}\};
\endmultline
$$
$P_{z_N}$ is pluripolar. Take a sequence
$A'\setminus(S_N\cup P_{z_N})\ni z'{}^\nu\too z'$. Then $z'{}^\nu\in
T_{z_N}$. Thus $\wdht f_{z_N}(z'{}^\nu)=f(z'{}^\nu,z_N)$. Hence,
by continuity, $\wdht f_{z_N}(z')=f(z',z_N)$\;
\ftnt{Here is the only place where the continuity of $f$ is used.}.

\halfskip

Moreover, $g\in\CO_s(Z)$. Put $V:=\MX(A',A_N;\wdht Y,D_N)\supset Z$.
Since the result is true for $N=2$ (without the continuity), we
get a holomorphic extension of $g$ to $\wdht V$. It remains to observe that
$\wdht V=\wdht X$; cf\. \cite{Jar-Pfl 2001c}, the proof of Step 3.
\qed\enddemo

\proclaim{Lemma 10} Let $D\subset\CC^p$, $G\subset\CC^q$ be pseudoconvex
domains, let $A\subset D$, $B\subset G$ be locally pluriregular, and let
$M\subset U$ be a relatively closed subset of an open neighborhood $U$ of
the cross $X:=\MX(A,B;D,G)$. Let $A'\subset A$, $B'\subset B$ be such that
$A\setminus A'$, $B\setminus B'$ are pluripolar and for any
$(a,b)\in A'\times B'$ the fibers $M_{(a,\cdot)}$, $M_{(\cdot,b)}$ are
pluripolar. Let $(D_j)_{j=1}^\infty$, $(G_j)_{j=1}^\infty$ be sequences of
pseudoconvex domains, $D_j\Subset D$, $G_j\Subset G$, with $D_j\nearrow D$,
$G_j\nearrow G$, such that $A'_j:=A'\cap D_j\neq\varnothing$,
$B'_j:=B'\cap G_j\neq\varnothing$, $j\in\NN$. We assume that for each
$j\in\NN$, $a\in A'_j$, and $b\in B'_j$, there exist:

$\bullet$ polydiscs $\Delta_a(r_{a,j})\subset D_j$,
$\Delta_b(s_{b,j})\subset G_j$ and

$\bullet$ relatively closed pluripolar sets
$S_{a,j}\subset\Delta_a(r_{a,j})\times G_j$,
$S^{b,j}\subset D_j\times\Delta_b(s_{b,j})$

\noindent such that:

$\bullet$ $(\Delta_a(r_{a,j})\times G_j)\cup
(D_j\times\Delta_b(s_{b,j}))\subset U\cap\wdht X$,

$\bullet$ $((A'\cap\Delta_a(r_{a,j}))\times G_j)\cap S_{a,j}\subset M$,
$(D_j\times(B'\cap\Delta_b(s_{b,j})))\cap S^{b,j}\subset M$,

$\bullet$ for any $f\in\CO_s(X\setminus M)$ there exist functions
$f_{a,j}\in\CO(\Delta_a(r_{a,j})\times G_j\setminus S_{a,j})$,
$f^{b,j}\in\CO(D_j\times\Delta_b(s_{b,j})\setminus S^{b,j})$ with

$f_{a,j}=f$ on $(A'\cap\Delta_a(r_{a,j}))\times G_j\setminus M$,

$f^{b,j}=f$ on $D_j\times(B'\cap\Delta_b(s_{b,j}))\setminus M$,

$\bullet$ $S_{a,j}$ is singular with respect to the family
$\{f_{a,j}\: f\in\CO_s(X\setminus M)\}$,
$S^{b,j}$ is singular with respect to the family
$\{f^{b,j}\: f\in\CO_s(X\setminus M)\}$.

\halfskip

Then there exists a relatively closed pluripolar set $\wdht M\subset\wdht X$
such that:

$\bullet$ $\wdht M\cap X'\subset M$, where $X':=\MX(A',B';D,G)$,

$\bullet$ for any $f\in\CO_s(X\setminus M)$ there exists exactly one
$\wdht f\in\CO(\wdht X\setminus\wdht M)$ with $\wdht f=f$ on $X'\setminus M$,

$\bullet$ the set $\wdht M$ is singular with respect to the family
$\{\wdht f\: f\in\CO_s(X\setminus M)\}$.
\endproclaim

\demo{Proof} Fix a $j\in\NN$. Put
$$
\align
\wdtl U_j:&=\bigcup_{a\in A'_j,\; b\in B'_j}
(\Delta_a(r_{a,j})\times G_j)\cup(D_j\times\Delta_b(s_{b,j})),\\
X_j:&=((A\cap D_j)\times G_j)\cup(D_j\times(B\cap G_j)),\\
X'_j:&=(A'_j\times G_j)\cup(D_j\times B'_j).
\endalign
$$
Note that $X'_j\subset\wdtl U_j$. Take an $f\in\CO_s(X\setminus M)$.
We like to glue the sets $(S_{a,j})_{a\in A'_j}$,
$(S^{b,j})_{b\in B'_j}$ and the functions $(f_{a,j})_{a\in A'_j}$,
$(f^{b,j})_{b\in B'_j}$ to obtain a global holomorphic function
$f_j:=\bigcup_{a\in A'_j,\; b\in B'_j}f_{a,j}\cup f^{b,j}$ on
$\wdtl U_j\setminus S_j$ where
$S_j:=\bigcup_{a\in A'_j,\; b\in B'_j}S_{a,j}\cup S^{b,j}$.

Let $a\in A'_j$, $b\in B'_j$. Observe that
$$
\align
f_{a,j}&=f \;\text{ on } (A'\cap\Delta_a(r_{a,j}))\times G_j\setminus M,\\
f^{b,j}&=f \;\text{ on } D_j\times(B'\cap\Delta_b(s_{b,j}))\setminus M.
\endalign
$$
Thus $f_{a,j}=f^{b,j}$ on the non-pluripolar set
$(A'\cap\Delta_a(r_{a,j}))\times (B'\cap\Delta_b(s_{b,j}))\setminus M$
(cf\. Lemma 8(b)). Hence
$$
f_{a,j}=f^{b,j}\;\text{ on }
\Delta_{a}(r_{a,j})\times\Delta_{b}(s_{b,j})\setminus(S_{a,j}\cup S^{b,j}).
$$
Using the minimality of $S_{a,j}$ and $S^{b,j}$, we conclude that
$$
S_{a,j}\cap(\Delta_{a}(r_{a,j})\times\Delta_{b}(s_{b,j}))=
S^{b,j}\cap(\Delta_{a}(r_{a,j})\times\Delta_{b}(s_{b,j})).
$$

Now let $a', a''\in A'_j$ be such that
$C:=\Delta_{a'}(r_{a',j})\cap\Delta_{a''}(r_{a'',j})\neq\varnothing$.
Fix a $b\in B'_j$. We know that $f_{a',j}=f^{b,j}=f_{a'',j}$ on
$C\times\Delta_b(r_{b,j})\setminus(S_{a',j}\cup S^{b,j}\cup S_{a'',j})$.
Hence, by the identity principle, we conclude that $f_{a',j}=f_{a'',j}$ on
$C\times G_j\setminus (S_{a',j}\cup S_{a'',j})$ and, moreover,
$$
S_{a',j}\cap(C\times G_j)=S_{a'',j}\cap(C\times G_j).
$$

The same argument works for $b', b''\in B'\cap G_j$.

\halfskip

Let $U_j$ be the connected component of $\wdtl U_j\cap\wdht X'_j$ with
$X'_j\subset U_j$. We have constructed a relatively closed pluripolar set
$S_j\subset U_j$ such that:

$S_j\cap X'_j\subset M$ and

for any $f\in\CO_s(X\setminus M)$ there exists (exactly one)
$f_j\in\CO(U_j\setminus S_j)$ with $f_j=f$ on $X'_j\setminus M$.

Recall that $X'_j\subset U_j\subset\wdht X'_j$. Hence the envelope of
holomorphy $\wdht U_j$ coincides with $\wdht X'_j$ (cf\. \cite{Jar-Pfl 2001b},
the proof of Step 4).

Applying the Chirka theorem (Theorem 6), we find a relatively closed
pluripolar set $\wdht M_j\subset\wdht X'_j$ such that:

$\wdht M_j\cap U_j\subset S_j$,

for any $f\in\CO_s(X\setminus M)$ there exists (exactly one) function
$\wdht f_j\in \CO(\wdht X'_j\setminus\wdht M_j)$ with
$\wdht f_j=f_j$ on $U_j\setminus S_j$ (in particular,
$\wdht f_j=f$ on $X'_j\setminus M$),

the set $\wdht M_j$ is singular with respect to the family
$\{\wdht f_j\: f\in\CO_s(X\setminus M)\}$.

\halfskip

Since $A\setminus A'$, $B\setminus B'$ are pluripolar, we get
$$
\align
\wdht X'_j&=\{(z,w)\in D_j\times G_j\: h_{A'\cap D_j, D_j}^\ast(z)+
h_{B'\cap G_j, G_j}^\ast(w)<1\}\\
&=\{(z,w)\in D_j\times G_j\: h_{A\cap D_j, D_j}^\ast(z)+
h_{B\cap G_j, G_j}^\ast(w)<1\}=\wdht X_j.
\endalign
$$ So, in fact, $\wdht f_j\in\CO(\wdht X_j\setminus\wdht M_j)$.
Observe that $\bigcup_{j=1}^\infty X_j=X$, $\wdht X_j\subset\wdht
X_{j+1}$, and $\bigcup_{j=1}^\infty\wdht X_j=\wdht X$. Using again
the minimality of the $\wdht M_j$'s (and gluing the $\wdht
f_j$'s), we get a relatively closed pluripolar set $\wdht
M\subset\wdht X$ which satisfies all the required conditions.
\qed\enddemo

\proclaim{Lemma 11}
Let $A\subset E^{n-1}$ be locally pluriregular, let $G\subset\CC$
be a domain with $E\Subset G$, let $X:=\MX(A,E;E^{n-1},G)$, and let
$U\subset E^{n-1}\times G$ be an open neighborhood of $X$. Let $M\subset U$
be a relatively closed set such that $M\cap E^n=\varnothing$ and for any
$a\in A$ the fiber $M_{(a,\cdot)}$ is polar. Then there exists a relatively
closed pluripolar set $\wdht M\subset\wdht X$ such that:

$\bullet$ $\wdht M\cap X\subset M$,

$\bullet$ for any $f\in\CO_s(X\setminus M)$ there exists exactly one
$\wdht f\in\CO(\wdht X\setminus\wdht M)$ with $\wdht f=f$ on $X\setminus M$,

$\bullet$ the set $\wdht M$ is singular with respect to the family
$\{\wdht f\: f\in\CO_s(X\setminus M)\}$.
\endproclaim

Notice that the above result is a special case of our Main Theorem with
$N=2$, $n_1=n-1$, $D_1=E^{n-1}$, $A_1=A$, $n_2=1$, $D_2=G$, $A_2=E$,
$\Sigma_1=\Sigma_2=\varnothing$.

\demo{Proof}
By Lemma 10, it suffices to show that for any $a_0\in A$ and for any domain
$G'\Subset G$ with $E\Subset G'$ there exist $r>0$ and a relatively closed
pluripolar set $S\subset\Delta_{a_0}(r)\times G'\subset U$ such that:

$\bullet$ $S\cap X\subset M$ and

$\bullet$ any function from $\CO_s(X\setminus M)$ extends
holomorphically to $\Delta_{a_0}(r)\times G'\setminus S$.

\halfskip

Fix $a_0$ and $G'$.
For $b\in G$, let $\rho=\rho_b>0$ be such that $\Delta_b(\rho)\Subset G$ and
$M_{(a_0,\cdot)}\cap\partial\Delta_b(\rho)=\varnothing$
(cf\. \cite{Arm-Gar 2001}, Th\. 7.3.9). Take $\rho^-=\rho_b^->0$,
$\rho^+=\rho_b^+>0$ such that $\rho^-<\rho<\rho^+$, $\Delta_b(\rho^+)\Subset
G$, and $M_{(a_0,\cdot)}\cap\overline P=\varnothing$, where
$$
P=P_b:=\{w\in\CC\: \rho^-<|w|<\rho^+\}.
$$
Let $\gamma\:[0,1]\too G\setminus M_{(a_0,\cdot)}$ be a curve such that
$\gamma(0)=0$, $\gamma(1)\in\partial\Delta_b(\rho)$. There exists an
$\eps=\eps_b>0$ such that
$$
\Delta_{a_0}(\eps)\times((\gamma([0,1])+\Delta_0(\eps))\cup P)\subset
U\setminus M.
$$
Put $V=V_b:=E\cup\big(\gamma([0,1])+\Delta_0(\eps)\big)\cup P$ and consider
the cross
$$
Y=Y_b:=\MX(A\cap\Delta_{a_0}(\eps),E;\Delta_{a_0}(\eps),V).
$$
Then $f\in\CO_s(Y)$ for any $f\in\CO_s(X\setminus M)$. Consequently, by
Theorem 3, any function from $\CO_s(X\setminus M)$ extends
holomorphically to $\wdht Y\supset\{a_0\}\times V$. Shrinking $\eps$ and $V$,
we may assume that any function $f\in\CO_s(X\setminus M)$ extends to a
function $\wdtl f=\wdtl f_b\in\CO(\Delta_{a_0}(\eps)\times W)$, where
$$
W=W_b:=\Delta_0(1-\eps)\cup(\gamma([0,1])+\Delta_0(\eps))\cup P.
$$

In particular, $\wdtl f$ is holomorphic in $\Delta_{a_0}(\eps)\times P$,
and therefore may be represented by the Hartogs--Laurent series
$$
\multline
\wdtl f(z,w)=\sum_{k=0}^\infty\wdtl f_k(z)(w-b)^k+
\sum_{k=1}^\infty\wdtl f_{-k}(z)(w-b)^{-k}=:\wdtl f^+(z,w)+\wdtl f^-(z,w),\\
(z,w)\in \Delta_{a_0}(\eps)\times P,
\endmultline
$$
where $\wdtl f^+\in\CO(\Delta_{a_0}(\eps)\times\Delta_b(\rho^+))$ and
$\wdtl f^-\in\CO(\Delta_{a_0}(\eps)\times(\CC\setminus
\overline\Delta_b(\rho^-)))$.
Recall that for any $a\in A\cap\Delta_{a_0}(\eps)$ the function
$\wdtl f(a,\cdot)$ extends holomorphically to $G\setminus M_{(a,\cdot)}$.
Consequently, for any $a\in A\cap\Delta_{a_0}(\eps)$ the function
$\wdtl f^-(a,\cdot)$ extends holomorphically to
$\CC\setminus(M_{(a,\cdot)}\cap\overline\Delta_b(\rho^-))$. Now, by Theorem 7,
there exists a relatively closed pluripolar set
$S=S_b\subset\Delta_{a_0}(\eps)\times\overline\Delta_b(\rho^-)$ such that:

$S\cap((A\cap\Delta_{a_0}(\eps))\times\overline\Delta_b(\rho^-)) \subset M$
and

any function $\wdtl f^-$ extends holomorphically to a function
$\overset\approx\to{f}^-\in\CO(\Delta_{a_0}(\eps)\times\CC\setminus S)$.

Since $\wdtl f=\wdtl f^++\wdtl f^-$, the function $\wdtl f$ extends holomorphically
to a function
$\wdht f=\wdht f_b\in\CO(\Delta_{a_0}(\eps)\times\Delta_b(\rho^+)\setminus S)$.
We may assume that the set $S$ is singular with respect to the
family $\{\wdht f\: f\in\CO_s(X\setminus M)\}$.

Using the identity principle and the minimality of the $S_b$'s, one can
easily show that for $b', b''\in G$, if
$B:=\Delta_{b'}(\rho_{b'}^+)\cap\Delta_{b''}(\rho_{b''}^+)\neq\varnothing$,
then
$$
S_{b'}\cap(\Delta_{a_0}(\eta)\times B)=S_{b''}\cap(\Delta_{a_0}(\eta)\times B),
\quad \wdht f_{b'}=\wdht f_{b''} \text{ on } \Delta_{a_0}(\eta)\times B,
$$
where $\eta:=\min\{\eps_{b'},\eps_{b''}\}$. Thus the functions
$\wdht f_{b'}$, $\wdht f_{b''}$ and sets $S_{b'}$, $S_{b''}$ may be glued
together.

Now, select $b_1,\dots,b_k\in G$ so that $G'\subset\bigcup_{j=1}^k
\Delta_{b_j}(\rho_{b_j}^+)$. Put
$$
r:=\min\{\eps_{b_j}\: j=1,\dots,k\}.
$$
Then $S:=(\Delta_{a_0}(r)\times G')\cap\bigcup_{j=1}^kS_{b_j}$ gives the
required relatively closed pluripolar subset of $\Delta_{a_0}(r)\times G'$
such that $S\cap X\subset M$ and for any $f\in\CO_s(X\setminus M)$,
the function $\wdht f:=\bigcup_{j=1}^k\wdht f_{b_j}$ extends holomorphically
$f$ to $\Delta_{a_0}(r)\times G'\setminus S$.
\qed\enddemo

\proclaim{Lemma 12}
Let $A\subset E^p$ be locally pluriregular, let $R>1$,
let
$$
X:=\MX(A,E^q;E^p,\Delta_0^q(R)),
$$
and let $U\subset E^p\times\Delta_0^q(R)$ be an open neighborhood of $X$. Let
$M\subset U$ be a relatively closed set such that $M\cap E^{p+q}=\varnothing$
and for any $a\in A$ the fiber $M_{(a,\cdot)}$ is pluripolar. Then there
exists a relatively closed pluripolar set $\wdht M\subset\wdht X$ such that:

$\bullet$ $\wdht M\cap X\subset M$,

$\bullet$ for any $f\in\CO_s(X\setminus M)$ there exists exactly one
$\wdht f\in\CO(\wdht X\setminus\wdht M)$ with $\wdht f=f$ on $X\setminus M$,

$\bullet$ the set $\wdht M$ is singular with respect to the family
$\{\wdht f\: f\in\CO_s(X\setminus M)\}$.
\endproclaim

Notice that the above result is a special case of our Main Theorem with
$N=2$, $n_1=p$, $D_1=E^p$, $A_1=A$, $n_2=q$, $D_2=\Delta_0^q(R)$, $A_2=E^q$,
$\Sigma_1=\Sigma_2=\varnothing$.

\demo{Proof} The case $q=1$ follows from Lemma 11. Thus assume that $q\geq2$.
By Lemma 10, it suffices to show that for any $a_0\in A$ and for any
$R'\in(1,R)$ there exist $r=r_{R'}>0$ and a relatively closed pluripolar set
$S=S_{R'}\subset\Delta_{a_0}(r)\times\Delta_0^q(R')\subset U$ such that

$\bullet$ $S\cap X\subset M$,

$\bullet$ any function from $\CO_s(X\setminus M)$ extends
holomorphically to $\Delta_{a_0}(r)\times\Delta_0^q(R')\setminus S$.

\halfskip

Fix an $a_0\in A$ and let $R'_0$ be the supremum of all $R'\in(0,R)$ such that
$r_{R'}$ and $S_{R'}$ exist. Note that $1\leq R'_0\leq R$. It suffices to show
that $R'_0=R$.

Suppose that $R'_0<R$. Fix $R'_0<R''<R$ and choose $R'\in(0,R'_0)$ such that
$\root{q}\of{{R'}^{q-1}R''}>R'_0$. Let $r:=r_{R'}$, $S:=S_{R'}$.

Write $w=(w',w_q)\in\CC^q=\CC^{q-1}\times\CC$. Let $C$ denote the set of all
$(a,b')\in(A\cap\Delta_{a_0}(r))\times\Delta_0^{q-1}(R')$ such that the fiber
$(M\cup S)_{(a,b',\cdot)}$ is polar. By Lemma 8(a,c), $C$ is pluriregular.
Now, by Lemma 11 applied to the cross
$$
Y_q:=\MX(C,\Delta_0(R');\Delta_{a_0}(r)\times\Delta_0^{q-1}(R'),\Delta_0(R))
$$
and the set $M_q:=M\cup S$, we conclude that there exists a closed pluripolar
set $S_q\subset\wdht Y_q$ such that $S_q\cap Y_q\subset M_q$ and any function
$f\in\CO_s(X\setminus M)$ extends holomorphically to $\wdht Y_q\setminus S_q$.
Using the product property of the relative extremal function
(cf\. \cite{Ngu-Sic 1991}), we get
$$
\align
\wdht Y_q&=\{(z,w',w_q)\in\Delta_{a_0}(r)\times\Delta_0^{q-1}(R')
\times\Delta_0(R)\:\\
&\hskip116pt h^\ast_{C,\Delta_{a_0}(r)\times\Delta_0^{q-1}(R')}(z,w')+
h^\ast_{\Delta_0(R'),\Delta_0(R)}(w_q)<1\}\\
&=\{(z,w',w_q)\in\Delta_{a_0}(r)\times\Delta_0^{q-1}(R')\times\Delta_0(R)\:\\
&\hskip40pt h^\ast_{(A\cap\Delta_{a_0}(r))\times\Delta_0^{q-1}(R'),
\Delta_{a_0}(r)\times\Delta_0^{q-1}(R')}(z,w')+
h^\ast_{\Delta_0(R'),\Delta_0(R)}(w_q)<1\}\\
&=\{(z,w',w_q)\in\Delta_{a_0}(r)\times\Delta_0^{q-1}(R')\times\Delta_0(R)\:\\
&\quad \max\{h^\ast_{A\cap\Delta_{a_0}(r),\Delta_{a_0}(r)}(z),\;
h^\ast_{\Delta_0^{q-1}(R'),\Delta_0^{q-1}(R')}(w')\}+
h^\ast_{\Delta_0(R'),\Delta_0(R)}(w_q)<1\}\\
&=\{(z,w',w_q)\in\Delta_{a_0}(r)\times\Delta_0^{q-1}(R')\times\Delta_0(R)\:\\
&\hskip140pt h^\ast_{A\cap\Delta_{a_0}(r),\Delta_{a_0}(r)}(z)+
h^\ast_{\Delta_0(R'),\Delta_0(R)}(w_q)<1\}.
\endalign
$$ Since $R''<R$, we find an $r_q\in(0,r]$ such that any function
$f\in\CO_s(X\setminus M)$ extends holomorphically to a function
$\wdtl f_q$ on
$\Delta_{a_0}(r_q)\times\Delta_0^{q-1}(R')\times\Delta_0(R'')\setminus
S_q$. We may assume that $S_q$ is singular with respect to the
family $\{\wdtl f_q: f\in\CO_s(X\setminus M)\}$.

\halfskip

Repeating the above argument for the coordinates $w_\nu$, $\nu=1,\dots,q-1$,
and gluing the obtained sets, we find an $r_0\in(0,r]$ and a relatively closed
pluripolar set $S_0:=\bigcup_{j=1}^qS_j$ such that any function
$f\in\CO_s(X\setminus M)$ extends holomorphically to a function
$\wdtl f_0:=\bigcup_{j=1}^q\wdtl f_j$ holomorphic in
$\Delta_{a_0}(r_0)\times\Omega\setminus S_0$, where
$$
\Omega:=\bigcup_{\nu=1}^q
\Delta_0^{j-1}(R')\times\Delta_0(R'')\times\Delta_0^{q-j}(R').
$$
Let $\wdht\Omega$ denote the envelope of holomorphy of $\Omega$. Applying the
Chirka theorem (Theorem 6), we find a relatively closed pluripolar subset
$\wdht S_0$ of $\Delta_{a_0}(r_0)\times\wdht\Omega$ such that any function
$f\in\CO_s(X\setminus M)$ extends to a function $\wdht f$ holomorphic on
$\Delta_{a_0}(r_0)\times\wdht\Omega\setminus\wdht S_0$. Let
$R''':=\root{q}\of{{R'}^{q-1}R''}$. Observe that
$\Delta_0(R''')\subset\wdht\Omega$. Recall that $R'''>R'_0$. We may assume
that $\wdht M$ is singular with respect to the family
$\{\wdht f: f\in\CO_s(X\setminus M)\}$. To get a contradiction it suffices to
show that $\wdht M\cap X\subset M$. We argue as in the proof of Lemma 11:

Take $(a,b)\in (A\cap\Delta_{a_0}(r_0))\times\Delta_0^q(R''')\setminus M$.
Since $M_{(a,\cdot)}$ is pluripolar, there exists a curve
$\gamma\:[0,1]\too\Delta_0(R''')\setminus M_{(a,\cdot)}$
such that $\gamma(0)=0$, $\gamma(1)=b$. Take an $\eps>0$ so small that
$$
\Delta_a(\eps)\times(\gamma([0,1])+\Delta_0^q(\eps))\subset
\Delta_{a_0}(r)\times\Delta_0^q(R''')\setminus M.
$$
Put $V_b:=E^q\cup(\gamma([0,1])+\Delta_0^q(\eps))$ and consider the cross
$$
Y:=\MX(A\cap\Delta_a(\eps),E^q;\Delta_a(\eps),V_b).
$$
Then $f\in\CO_s(Y)$ for any $f\in\CO_s(X\setminus M)$. Consequently, by
Theorem 3,
$\wdht Y\subset\Delta_{a_0}(r)\times\Delta_0^q(R''')\setminus\wdht M$, which
implies that $\wdht M_{(a,\cdot)}\cap\Delta_0^q(R''')\subset M_{(a,\cdot)}$.
\qed\enddemo

\halfskip

\noindent{\bf 3. Proof of the Main Theorem for $\bold{N=2}$.}
To simplify notation put:
$p:=n_1$, $D:=D_1$, $A:=A_1$, $A':=A\setminus\Sigma_2$,
$q:=n_2$, $G:=D_2$, $B:=A_2$, $B':=B\setminus\Sigma_1$.

It suffices to verify the assumptions of Lemma 10.
Let $(D_j)_{j=1}^\infty$, $(G_j)_{j=1}^\infty$ be approximation sequences:
$D_j\Subset D_{j+1}\Subset D$, $G_j\Subset G_{j+1}\Subset G$, $D_j\nearrow D$, $G_j\nearrow G$,
$A'\cap D_j\neq\varnothing$, and $B'\cap G_j\neq\varnothing$, $j\in\NN$.

Fix a $j\in\NN$, $a\in A'\cap D_j$ and let $\Omega_j$ be the set of all
$b\in G_{j+1}$ such that there exist a polydisc
$\Delta_{(a,b)}(r_b)\subset D_j\times G_{j+1}$ and a relatively
closed pluripolar set $S_b\subset\Delta_{(a,b)}(r_b)$ such that:

$S_b\cap((A'\cap\Delta_a(r_b))\times\Delta_b(r_b))\subset M$,

any function $f\in\CO_s(X\setminus M)$ extends to a function
$\wdtl f_b\in\CO(\Delta_{(a,b)}(r_b)\setminus S_b)$ with $\wdtl f_b=f$ on
$(A'\cap\Delta_a(r_b))\times\Delta_b(r_b)\setminus M$,

and $S_b$ is singular with respect to the family
$\{\wdtl f_b\: f\in\CO_s(X\setminus M)\}$.

It is clear that $\Omega_j$ is open. Observe that $\Omega_j\neq\varnothing$.
Indeed, since $B\cap G_j\setminus M_{(a,\cdot)}\neq\varnothing$,
we find a point $b\in B\cap G_j\setminus M_{(a,\cdot)}$. Therefore there is a polydisc
$\Delta_{(a,b)}(r)\subset D_j\times G_j\setminus M$. Put
$$
Y:=\MX(A\cap\Delta_a(r),B\cap\Delta_b(r);\Delta_a(r),\Delta_b(r)).
$$
By Theorem 3, we find
an $r_b\in(0,r)$ such that any function $f\in\CO_s(X\setminus M)$ extends
to $\wdtl f_b\in\CO(\Delta_{(a,b)}(r_b))$ with
$\wdtl f_b=f$ on $\Delta_{(a,b)}(r_b)\cap Y\supset
(A\cap\Delta_a(r_b))\times\Delta_b(r_b)$. Consequently, $b\in\Omega_j$.

Moreover, $\Omega_j$ is relatively closed in $G_{j+1}$.
Indeed, let $c$ be an accumulation point of $\Omega_j$ in $G_{j+1}$ and
let $\Delta_c(3R)\subset G_{j+1}$. Take a point $b\in\Omega_j\cap\Delta_c(R)
\setminus M_{(a,\cdot)}$ and let $r\in(0,r_b]$, $r<2R$,
be such that $\Delta_{(a,b)}(r)\cap M=\varnothing$.
Observe that $\wdtl f_b\in\CO(\Delta_{(a,b)}(r))$ and
$\wdtl f_b(z,\cdot)=f(z,\cdot)\in\CO(\Delta_b(2R)\setminus M_{(z.\cdot)})$
for any $z\in A'\cap\Delta_a(r)$. Hence, by Lemma 12 (with $R':=R$), there
exists a relatively closed pluripolar set
$S\subset\Delta_a(\rho')\times\Delta_b(R)$ with $\rho'\in(0,r)$ such that any
$f$ has an extension
$\wdht{\wdtl f_b}\in\CO(\Delta_a(\rho')\times\Delta_b(R)\setminus S)$
Take an $r_c>0$ so small that
$\Delta_{(a,c)}(r_c)\subset\Delta_a(\rho')\times\Delta_b(R)$ and put
$S_c:=S\cap\Delta_{(a,c)}(r_c)$,
$\wdtl f_c:=\wdht{\wdtl f_b}$ on $\Delta_{(a,c)}(r_c)\setminus S_c$.
Obviously $\wdtl f_c=\wdht{\wdtl f_b}=f$ on $(A'\cap\Delta_a(r_c))
\times\Delta_c(r_c)\setminus M$. Hence $c\in\Omega_j$.

Thus $\Omega_j=G_{j+1}$. There exists a finite set $T\subset\overline G_j$
such that
$$
\overline G_j\subset\bigcup_{b\in T}\Delta_b(r_b).
$$
Define $r_{a,j}:=\min\{r_b\: b\in T\}$. Take $b', b''\in T$ with
$C:=\Delta_{b'}(r_{b'})\cap\Delta_{b''}(r_{b''})\neq\varnothing$.
Then $\wdtl f_{b'}=f=\wdtl f_{b''}$ on $(A'\cap\Delta_a(r_{a,j}))
\times(\Delta_{b'}(r_{b'})\cap\Delta_{b''}(r_{b''}))\setminus M$.
Consequently,  $\wdtl f_{b'}=\wdtl f_{b''}$ on
$\Delta_a(r_{a,j})\times C\setminus(S_{b'}\cup S_{b''})$.
In particular, using the minimality of the sets $S_{b'}$ and $S_{b''}$, we
conclude that they coincide on $\Delta_a(r_{a,j})\times C$ and that
the functions $f_{b'}$ and $f_{b''}$ glue together.
Thus we get a relatively closed pluripolar set $S_{a,j}\subset
\Delta_a(r_{a,j})\times G_j$ such that
$S_{a,j}\cap((A'\cap\Delta_a(r_{a,j}))\times G_j)\subset M$ and
any function $f\in\CO_s(X\setminus M)$ extends holomorphically to
an $f_{a,j}\in\CO(\Delta_a(r_{a,j})\times G_j\setminus S_{a,j})$ with
$f_{a,j}=f$ on $(A'\cap\Delta_a(r_{a,j}))\times G_j\setminus M$.

Changing the role of $z$ and $w$, we get $S^{b,j}$ and
$f^{b,j}$, $b\in B'\cap G_j$.
\qed

\halfskip

The above proof of the Main Theorem for $N=2$ shows that the following
generalization of Lemma 12 is true.

\proclaim{Theorem 13}
Let $D\subset\CC^p$, $G\subset\CC^q$ be pseudoconvex
domains, let $A\subset D$ be locally pluriregular, let $B\subset G$ be open
and non-empty, and let $M\subset U$ be a relatively closed subset of an open
neighborhood $U$ of the cross $X:=\MX(A,B;D,G)$ such that
$M\cap(D\times B)=\varnothing$ and for any $a\in A$ the fiber $M_{(a,\cdot)}$ is pluripolar.
Then there exists a relatively closed pluripolar set $\wdht M\subset\wdht X$
such that:

$\bullet$ $\wdht M\cap X\subset M$,

$\bullet$ for any $f\in\CO_s(X\setminus M)$ there exists exactly one
$\wdht f\in\CO(\wdht X\setminus\wdht M)$ with $\wdht f=f$ on
$X\setminus M$,

$\bullet$ the set $\wdht M$ is singular with respect to the family
$\{\wdht f\: f\in\CO_s(X\setminus M)\}$.
\endproclaim

Observe that if $G=\CC^q$, then $\wdht X=D\times\CC^q$. Consequently,
Theorem 13 generalizes also Theorem 7.

\demo{Proof} We apply Lemma 10 (as in the proof of the Main Theorem for $N=2$).
The functions $f_{a,j}$ are constructed exactly as in that proof (with $A'=A$).
The functions $f^{b,j}$ are simply given as
$f^{b,j}:=f|_{D_j\times\Delta_b(s_{b,j})}$
with $\Delta_b(s_{b,j})\subset B\cap D_j$ ($S^{b,j}:=\varnothing$).
\qed\enddemo


\halfskip
\noindent{\bf 4. Proof of the Main Theorem.}
First observe that, by Lemma 8(b), the set $X'\setminus M$ is not pluripolar.
Consequently, the function $\wdht f$ is uniquely determined.

We proceed by induction on $N$. The case $N=2$ is proved.

Let $D_{j,k}\nearrow D_j$, $D_{j,k}\Subset D_{j,k+1}\Subset D_j$, where
$D_{j,k}$ are pseudoconvex domains with
$A_{j,k}:=A_j\cap D_{j,k}\neq\varnothing$, $j=1,\dots,N$. Put
$$
\gather
X_k:=\MX(A_{1,k},\dots,A_{N,k};D_{1,k},\dots,D_{N,k})\subset X,\\
\Sigma_{j,k}:=(A_{1,k}\tdots A_{j-1,k}\times A_{j+1,k}\tdots A_{N,k})
\cap\Sigma_j,\quad j=1,\dots,N,\\
X'_k:=\MT(A_{1,k},\dots,A_{N,k};D_{1,k},\dots,D_{N,k};\Sigma_{1,k},\dots,
\Sigma_{N,k})\subset X_k.
\endgather
$$

\halfskip

It suffices to show that for each $k\in\NN$ the following condition
(*) holds.

\halfskip

\noindent(*) \quad There exists a domain $U_k$,
$X'_k\subset U_k\subset\wdht X_k$ and a relatively closed pluripolar
set $M_k\subset U_k$, such that:

$M_k\cap X'_k\subset M$ and

for any $f\in\CO_s(X\setminus M)$ there exists an
$\wdtl f_k\in\CO(U_k\setminus M_k)$ with $\wdtl f_k=f$ on $X'_k\setminus M$.

\halfskip

Indeed, fix a $k\in\NN$ and observe that, by Lemma 9,
$\wdht X_k$ is the envelope of holomorphy of $U_k$. Hence, in virtue of the
Chirka theorem (Theorem 6), there exists a relatively closed pluripolar
set $\wdht M_k$ of $\wdht X_k$, $\wdht M_k\cap U_k\subset M_k$, such that
$\wdht X_k\setminus\wdht M_k$ is the envelope of holomorphy of
$U_k\setminus M_k$. In particular, for each $f\in\CO_s(X\setminus M)$ there
exists an $\wdht f_k\in\CO(\wdht X_k\setminus\wdht M_k)$ with
$\wdht f_k|_{U_k\setminus M_k}=\wdtl f_k$.
We may assume that $\wdht M_k$ is singular with respect to the family
$\{\wdht f_k\: f\in\CO_s(X\setminus M)\}$.

\halfskip

In particular, $\wdht M_{k+1}\cap\wdht X_k=\wdht M_k$. Consequently:

$\wdht M:=\bigcup_{k=1}^\infty\wdht M_k$ is a relatively closed pluripolar
subset of $\wdht X$ with $\wdht M\cap X'\subset M$,

for each $f\in\CO_s(X\setminus M)$, the function
$\wdht f:=\bigcup_{k=1}^\infty\wdht f_k$ is holomorphic on
$\wdht X\setminus\wdht M$ with $\wdht f=f$ on $X'\setminus M$, and

$\wdht M$ is singular with respect to the family
$\{\wdht f\: f\in\CO_s(X\setminus M)\}$.

\halfskip

It remains to prove (*). Fix a $k\in\NN$. For any
$$
a=(a_1,\dots,a_N)\in A_{1,k}\tdots A_{N,k}\setminus M
$$
let $\tau=\tau_k(a)$ be such that
$\Delta_a(\tau)\subset D_{1,k}\tdots D_{N,k}\setminus M$. Consider
the $N$--fold cross
$$
Y_a:=\MX(A_1\cap\Delta_{a_1}(\tau),\dots,A_N\cap\Delta_{a_N}(\tau);
\Delta_{a_1}(\tau),\dots,\Delta_{a_N}(\tau)).
$$
Observe that any function from $\CO_s(X\setminus M)$ belongs to
$\CO_s(Y_a)$. Consequently, by Theorem 3, any function
from $\CO_s(X\setminus M)$ extends holomorphically to $\wdht Y_a$.
Let $\rho=\rho_k(a)\in(0,\tau]$ be such that $\Delta_a(\rho)\subset\wdht Y_a$.

If $N\geq4$, then we additionally define $(N-2)$--fold crosses
$$
\align
Y_{k,\mu,\nu}:=\MX(&A_{1,k},\dots,A_{\mu-1,k},A_{\mu+1,k},\dots,
A_{\nu-1,k},A_{\nu+1,k},\dots,A_{N,k};\\
&D_{1,k},\dots,D_{\mu-1,k},D_{\mu+1,k},\dots,
D_{\nu-1,k},D_{\nu+1,k},\dots,D_{N,k}),\\
&\hskip200pt 1\leq\mu<\nu\leq N,
\endalign
$$
and we assume that $\rho$ is so small that
$$
\Delta_{(a_1,\dots,a_{\mu-1},a_{\mu+1},\dots, a_{\nu-1},a_{\nu+1},\dots,a_N)}
(\rho)\subset\wdht Y_{k,\mu,\nu},\quad 1\leq\mu<\nu\leq N.
$$

For $j\in\{1,\dots,N\}$, define the $2$--fold crosses: $$
\multline Z'_{k,a,j}:=\Big\{(z',z_j,z'')\in
((A_1\cap\Delta_{a_1}(\rho))
\tdots(A_{j-1}\cap\Delta_{a_{j-1}}(\rho)))\times D_{j,k+1}\\
\times((A_{j+1}\cap\Delta_{a_{j+1}}(\rho))\tdots(A_N\cap\Delta_{a_N}(\rho)))\:
(z',z'')\notin\Sigma_j\Big\}\cup\Delta_a(\rho),
\endmultline
$$
$$
\multline
Z_{k,a,j}:=
\Big((A_1\cap\Delta_{a_1}(\rho))\tdots(A_{j-1}\cap\Delta_{a_{j-1}}(\rho))
\times D_{j,k+1}\\
\times(A_{j+1}\cap\Delta_{a_{j+1}}(\rho))\tdots(A_N\cap\Delta_{a_N}
(\rho)\Big)\cup\Delta_a(\rho).
\endmultline
$$
Now, we apply Theorem 13 to the $2$--fold cross $Z'_{k,a,j}$
and the set $M$. We find a relatively closed
pluripolar set $S_{k,a,j}\subset\wdht Z'_{k,a,j}=\wdht Z_{k,a,j}$ such that:

$S_{k,a,j}\cap Z'_{k,a,j}\subset M$,

for any function $f\in\CO_s(X\setminus M)$ there exists an
$\wdtl f_{k,a,j}\in\CO(\wdht Z_{k,a,j}\setminus S_{k,a,j})$ such that
$\wdtl f_{k,a,j}=f$ on $Z'_{k,a,j}\setminus M$,

$S_{k,a,j}$ is singular with respect to the space
$\{\wdtl f_{k,a,j}\: f\in\CO_s(X\setminus M)\}$.

Observe that $\{(a_1,\dots,a_{j-1})\}\times\overline D_{j,k}\times
\{(a_{j+1},\dots,a_N)\}\Subset\wdht Z_{k,a,j}$. Consequently, we
find $r=r_k(a)\in(0,\rho]$ such that
$$
V_{k,a,j}:=\Delta_{(a_1,\dots,a_{j-1})}(r)\times D_{j,k}\times
\Delta_{(a_{j+1},\dots,a_N)}(r)\subset\wdht Z_{k,a,j},\quad j=1,\dots,N.
$$
Let
$$
V_k:=\bigcup_{\Sb a\in A_{1,k}\tdots A_{N,k}\setminus M\\ j\in\{1,\dots,N\}
\endSb}V_{k,a,j}.
$$
Note that $X'_k\subset V_k$. Let $U_k$ be the connected component of
$V_k\cap\wdht X_k$ that contains $X_k$.

\halfskip

It remains to glue the sets $S_{k,a,j}$ and functions $\wdtl f_{k,a,j}$. Then
$$
S_k:=\bigcup_{\Sb a\in A_{1,k}\tdots A_{N,k}\setminus M\\ j\in\{1,\dots,N\}
\endSb}S_{k,a,j}\cap U_k,\quad
\wdtl f_k:=\bigcup_{\Sb a\in A_{1,k}\tdots A_{N,k}\setminus M\\
j\in\{1,\dots,N\}\endSb}\wdtl f_{k,a,j}|_{V_{k,a,j}\cap U_k\setminus S_k}
$$
will satisfy (*).

\halfskip

To check that the gluing process is possible, let
$a, b\in A_{1,k}\tdots A_{N,k}\setminus M$, $i,j\in\{1,\dots,N\}$ be such that
$V_{k,a,i}\cap V_{k,b,j}\neq\varnothing$. We have the following two cases:

\halfskip

(a) $i\neq j$: We may assume that $i=N-1$, $j=N$. Write
$w=(w',w'')\in\CC^{n_1+\dots+n_{N-2}}\times\CC^{n_{N-1}+n_N}$. Observe that
$$
V_{k,a,N-1}\cap V_{k,b,N}
=\Big(\Delta_{a'}(r_k(a))\cap\Delta_{b'}(r_k(b))\Big)\times
\Delta_{b_{N-1}}(r_k(b))\times\Delta_{a_N}(r_k(a)).
$$

We consider the following three subcases:

\halfskip

$N=2$ (cf\. the proof of Lemma 10): Then $V_{k,a,1}\cap V_{k,b,2}
=\Delta_{b_1}(r_k(b))\times\Delta_{a_2}(r_k(a))$. We know that
$\wdtl f_{k,a,1}=\wdtl f_{k,b,2}$ on the non-pluripolar set
$$
(A_1\cap\Delta_{b_1}(r_k(b))\setminus\Sigma_2)\times
(A_2\cap\Delta_{a_2}(r_k(a))\setminus\Sigma_1)\setminus M;
$$
cf\. Lemma 8(b). Hence, by the identity principle,
$\wdtl f_{k,a,1}=\wdtl f_{k,b,2}$ on
$V_{k,a,1}\cap V_{k,b,2}\setminus (S_{k,a,1}\cup S_{k,b,2})$.
Consequently, the sets $S_{k,a,1}$, $S_{k,b,2}$ and the functions
$\wdtl f_{k,a,1}$, $\wdtl f_{k,b,2}$ glue together.

\halfskip

$N=3$: Then $V_{k,a,2}\cap V_{k,b,3}=
(\Delta_{a_1}(r_k(a))\cap\Delta_{b_1}(r_k(b))\times
\Delta_{b_2}(r_k(b))\times\Delta_{a_3}(r_k(a))$. Let
$$
C'':=(A_2\cap\Delta_{b_2}(r_k(b)))\times
(A_3\cap\Delta_{a_3}(r_k(a)))\setminus\Sigma_1.
$$
Recall that for any $c''\in C''$ the fiber $M_{(\cdot,c'')}$ is pluripolar.
We have $\wdtl f_{k,a,2}(\cdot,c'')=f(\cdot,c'')=\wdtl f_{k,b,3}(\cdot,c'')$
on $\Delta_{a_1}(r_k(a))\cap\Delta_{b_1}(r_k(b))\setminus M_{(\cdot,c'')}$.

Now, let $C'$ denote the set of all
$c'\in\Delta_{a_1}(r_k(a))\cap\Delta_{b_1}(r_k(b))$ such that the fiber
$(S_{k,a,2}\cup S_{k,b,3})_{(c',\cdot)}$ is pluripolar. Recall that the
complement of $C'$ is pluripolar (Lemma 8(a)). If $c'\in C'$,
then $\wdtl f_{k,a,2}(c',\cdot)=\wdtl f_{k,b,3}(c',\cdot)$
on $C''\setminus (S_{k,a,2}\cup S_{k,b,3})_{(c',\cdot)}$.
Consequently, by the identity principle,
$\wdtl f_{k,a,2}(c',\cdot)=\wdtl f_{k,b,3}(c',\cdot)$
on $\Delta_{b_2}(r_k(b))\times\Delta_{a_3}(r_k(a))\setminus
(S_{k,a,2}\cup S_{k,b,3})_{(c',\cdot)}$, $c'\in C'$.
Finally, $\wdtl f_{k,a,2}=\wdtl f_{k,b,3}$ on
$V_{k,a,2}\cap V_{k,b,3}\setminus(S_{k,a,2}\cup S_{k,b,3})$.
Consequently, the sets $S_{k,a,2}$, $S_{k,b,3}$ and the functions
$\wdtl f_{k,a,2}$, $\wdtl f_{k,b,3}$ glue together.

\halfskip

If $N\in\{2,3\}$, then we jump directly to (b) and we conclude that the Main
Theorem is true for $N\in\{2,3\}$.

\halfskip

$N\geq4$: Here is the only place where the induction over $N$ is used. We
assume that the Main Theorem is true for $N-1\geq3$.

Let
$$
\multline
C'':=\{c''\in(A_{N-1}\cap\Delta_{b_{N-1}}(r_k(b)))\times
(A_N\cap\Delta_{a_N}(r_k(a)))\:\\ (\Sigma_s)_{(\cdot,c'')}
\text{ is pluripolar},\;s=1,\dots,N-2\};
\endmultline
$$
note that, by Lemma 8(a), $C''$ is non-pluripolar. For any $c''\in C''$ the
function $f_{c''}:= f(\cdot,c'')$ is separately holomorphic on
$Y_{k,N-1,N}\setminus M_{(\cdot,c'')}$. Moreover, the set $M_{(\cdot,c'')}$
satisfies all the assumptions of the Main Theorem. Indeed,
$$
\multline
\Sigma_s(A_{1,k},\dots,A_{N-2,k};M_{(\cdot,c'')})=
(\Sigma_s(A_{1,k},\dots,A_{N,k};M))_{(\cdot,c'')}\subset
(\Sigma_s)_{(\cdot,c'')},\\
\quad s=1,\dots,N-2.
\endmultline
$$

By the inductive assumption, the
function $f_{c''}$ extends to a function
$$
\wdht f_{c''}\in\CO(\wdht Y_{k,N-1,N}\setminus\wdht M(c'')),
$$
where $\wdht M(c'')$ is
relatively closed pluripolar subset of $\wdht Y_{k,N-1,N}$ such that
$\wdht M(c'')\cap Y'_{k,N-1,N} \subset M_{(\cdot,c'')}$. Recall that
$$
\Delta_{a'}(r_k(a))\cup\Delta_{b'}(r_k(b))\subset\wdht Y_{k,N-1,N}.
$$
Since $\wdtl f_{k,a,N-1}(\cdot,c'')=f_{c''}$ on
$\Delta_{a'}(r_k(a))\cap Y'_{k,N-1,N}\setminus M_{(\cdot,c'')}$ and
$\wdtl f_{k,b,N}(\cdot,c'')=f_{c''}$ on
$\Delta_{b'}(r_k(b))\cap Y'_{k,N-1,N}\setminus M_{(\cdot,c'')}$, we conclude
that $\wdtl f_{k,a,N-1}(\cdot,c'')=\wdht f_{c''}=\wdtl f_{k,b,N}(\cdot,c'')$
on $\Delta_{a'}(r_k(a))\cap\Delta_{b'}(r_k(b))\setminus M_{(\cdot,c'')}$.

Let $c'\in\Delta_{a'}(r_k(a))\cap\Delta_{b'}(r_k(b))$ be such that the fiber
$(S_{k,a,N-1}\cup S_{k,b,N})_{(c',\cdot)}$ is pluripolar.
Then $\wdtl f_{k,a,N-1}(c',\cdot)=\wdtl f_{k,b,N}(c',\cdot)$
on $C''\setminus (S_{k,a,N-1}\cup S_{k,b,N})_{(\cdot, c')}$.
Consequently, by the identity principle,
$\wdtl f_{k,a,N-1}(c',\cdot)=\wdtl f_{k,b,N}(c',\cdot)$ on
$(\Delta_{b_{N-1}}(r_k(b))\times\Delta_{a_N}(r_k(a)))\setminus
(S_{k,a,N-1}\cup S_{k,b,N})_{c'}$ and,
finally, $\wdtl f_{k,a,N-1}=\wdtl f_{k,b,N}$ on
$(V_{k,a,N-1}\cap V_{k,b,N})\setminus(S_{k,a,N-1}\cup S_{k,b,N})$.
Consequently, the sets $S_{k,a,N-1}$, $S_{k,b,N}$ and the functions
$\wdtl f_{k,a,N-1}$, $\wdtl f_{k,b,N}$ glue together.

\halfskip

(b) $i=j$: We may assume that $i=j=N$. Observe that
$$
V_{k,a,N}\cap V_{k,b,N}=\Big(\Delta_{(a_1,\dots,a_{N-1})}(r_k(a))\cap
\Delta_{(b_1,\dots,b_{N-1})}(r_k(b))\Big)\times D_{N,k}.
$$
By (a) we know that:
$$
\alignat2
\wdtl f_{k,a,N}&=\wdtl f_{k,a,N-1} &&\text{ on }
V_{k,a,N}\cap V_{k,a,N-1}\setminus(S_{k,a,N}\cup S_{k,a,N-1}),\\
\wdtl f_{k,a,N-1}&=\wdtl f_{k,b,N} &&\text{ on }
V_{k,a,N-1}\cap V_{k,b,N}\setminus (S_{k,a,N-1}\cup S_{k,b,N}).
\endalignat
$$
Hence (we write $w=(w',w_N)\in\CC^{n_1+\dots+n_{N-1}}\times\CC^{n_N}$):
$$
\multline
\wdtl f_{k,a,N}=\wdtl f_{k,b,N} \text{ on }
V_{k,a,N}\cap V_{k,a,N-1}\cap V_{k,b,N}\setminus
(S_{k,a,N-1}\cup S_{k,a,N}\cup S_{k,b,N})\\
=\Big(\Delta_{a'}(r_k(a))\cap
\Delta_{b'}(r_k(b))\Big)\times\Delta_{a_N}(r_k(a))
\setminus (S_{k,a,N-1}\cup S_{k,a,N}\cup S_{k,b,N}),
\endmultline
$$
and finally, by the identity principle,
$$
\wdtl f_{k,a,N}=\wdtl f_{k,b,N} \text{ on }
V_{k,a,N}\cap V_{k,b,N}\setminus (S_{k,a,N}\cup S_{k,b,N}).
$$
Consequently, the sets $S_{k,a,N}$, $S_{k,b,N}$ and the functions
$\wdtl f_{k,a,N}$, $\wdtl f_{k,b,N}$ glue together.

\halfskip

The proof of the Main Theorem is completed.
\qed

\Refs

{\redefine\bf{\rm}
\widestnumber{\key}{XXXXXXXXX}
\ref
\key Ale-Zer 2001
\by O\. Alehyane, A\. Zeriahi
\paper Une nouvelle version du th\'eor\`eme d'extension de Hartogs pour les
applications s\'epar\'ement holomorphes entre espaces analytiques
\jour Ann\. Polon\. Math\.
\vol 76
\yr 2001
\pages 245--278
\endref
\ref
\key Arm-Gar 2001
\by D\. H\. Armitage, S\. J\. Gardiner
\book Classical Potential Theory
\publ Springer Verlag
\yr 2001
\endref
\ref
\key Chi 1993
\by E\. M\. Chirka
\paper The extension of pluripolar singularity sets
\jour Proc\. Steklov Inst\. Math\.
\vol 200
\yr 1993
\pages 369--373
\endref
\ref
\key Chi-Sad 1988
\by E\. M\. Chirka, A\. Sadullaev
\paper On continuation of functions with polar singularities
\jour Math\. USSR-Sb\.
\vol 60
\yr 1988
\pages 377--384
\endref
\ref
\key Jar-Pfl 2000
\by  M\. Jarnicki, P\. Pflug
\book  Extension of Holomorphic Functions
\publ  de Gruyter Expositions in Mathematics 34, Walter de Gruyter
\yr 2000
\endref
\ref
\key Jar-Pfl 2001a
\by  M\. Jarnicki, P\. Pflug
\paper Cross theorem
\jour Ann\. Polon\. Math\.
\vol 77
\yr 2001
\pages 295--302
\endref
\ref
\key Jar-Pfl 2001b
\by  M\. Jarnicki, P\. Pflug
\paper An extension theorem for separately holomorphic functions with
singularities
\jour IMUJ Preprint 2001/27
\yr 2001
\endref
\ref
\key Jar-Pfl 2001c
\by  M\. Jarnicki, P\. Pflug
\paper An extension theorem for separately holomorphic functions with
analytic singularities
\toappear
\endref
\ref
\key Kaz 1988
\by M\. V\. Kazaryan
\paper On holomorphic continuation of functions with pluripolar singularities
\jour Dokl\. Akad\. Nauk Arm\. SSR
\yr 1988
\vol 87
\pages 106--110 (in Russian)
\endref
\ref
\key Ngu 1997
\by Nguyen Thanh Van
\paper Separate analyticity and related subjects
\jour Vietnam J\. Math\.
\vol 25
\yr 1997
\pages 81--90
\endref
\ref
\key Ngu-Sic 1991
\by Nguyen Thanh Van, J\. Siciak
Fonctions plurisousharmoniques extr\'emales et syst\`emes doublement
orthogonaux de fonctions analytiques
\jour Bull\. Sc\. Math\.
\yr 1991
\vol 115
\pages 235--244
\endref
\ref
\key Ngu-Zer 1991
\by Nguyen Thanh Van, A\. Zeriahi
\paper  Une extension du th\'eor\`eme de
Hartogs sur les fonctions s\'epar\'ement analytiques
\inbook Analyse Complexe Multivariables, R\'ecents D\`evelopements,
A\. Meril (ed.), EditEl, Rende
\yr 1991
\pages 183--194
\endref
\ref
\key Ngu-Zer 1995
\by Nguyen Thanh Van, A\. Zeriahi
\paper Syst\`emes doublement othogonaux de fonctions holomorphes et applications
\jour Banach Center Publ\.
\vol 31
\yr 1995
\pages 281--297
\endref
\ref
\key \"Okt 1998
\by O\. \"Oktem
\paper Extension of separately analytic functions and applications to
range characterization of exponential Radon transform
\jour Ann\. Polon\. Math\.
\vol 70
\yr 1998
\pages 195--213
\endref
\ref
\key \"Okt 1999
\by O\. \"Oktem
\paper Extending separately analytic functions in $\CC^{n+m}$ with
singularities
\inbook Extension of separately analytic functions and applications to
mathematical tomography (Thesis)
\publ Dep\. Math\. Stockholm Univ\.
\yr 1999
\endref
\ref
\key Shi 1989
\by B\. Shiffman
\paper On separate analyticity and Hartogs theorem
\jour Indiana Univ\. Math\. J\.
\vol 38
\yr 1989
\pages 943--957
\endref
\ref
\key Sic 1969
\by J\. Siciak
\paper Separately analytic functions and envelopes
of holomorphy of some lower dimensional subsets of $\CC^n$
\jour Ann\. Polon\. Math\.
\vol 22
\yr 1969--1970
\pages 147--171
\endref
\ref
\key Sic 1981
\by J\. Siciak
\paper Extremal plurisubharmonic functions in $\CC^N$
\jour Ann\. Polon\. Math\.
\vol 39
\yr 1981
\pages 175--211
\endref
\ref \key Sic 2000
\by J\. Siciak
\paper Holomorphic functions with singularities on algebraic sets
\jour IMUJ Pre\-print 2000/21 (to appear in Univ\. Iag\. Acta Math\. 75 (2001))
\endref
\ref
\key Zah 1976
\by V\. P\. Zahariuta
\paper Separately analytic functions, generalizations of
Hartogs theorem, and envelopes of holomorphy
\jour Math\. USSR-Sb\.
\vol 30
\yr 1976
\pages 51--67
\endref
}
\endRefs

\enddocument